\documentclass[11pt]{amsart}
\usepackage{amscd,amssymb,subfigure}
\usepackage[arrow,matrix,graph]{xy}

\newtheorem{theorem}{Theorem}[section]
\newtheorem{corollary}[theorem]{Corollary}
\newtheorem{lemma}[theorem]{Lemma}
\newtheorem{prop}[theorem]{Proposition}

\theoremstyle{definition}

\newtheorem{example}[theorem]{Example}
\newtheorem{remark}[theorem]{Remark}
\newtheorem*{ack}{Acknowledgments}
\newtheorem*{claim}{Claim}
\newtheorem*{addproof}{Added in Proof}

\numberwithin{equation}{section}

\newcommand{\N}{\mathbb{N}}
\newcommand{\Z}{\mathbb{Z}}
\newcommand{\Q}{\mathbb{Q}}
\newcommand{\R}{\mathbb{R}}
\newcommand{\C}{\mathbb{C}}

\renewcommand{\k}{\Bbbk}
\newcommand{\HH}{{\mathfrak H}}

\newcommand{\RR}{{\mathcal R}^1}
\newcommand{\LL}{{\mathbf L}}

\newcommand{\G}{\Gamma}
\newcommand{\E}{\mathsf{E}}
\newcommand{\V}{\mathsf{V}}

\newcommand{\W}{\mathsf{W}}

\DeclareMathOperator{\rank}{rank}
\DeclareMathOperator{\gr}{gr}
\DeclareMathOperator{\im}{im}

\DeclareMathOperator{\ch}{char}
\DeclareMathOperator{\Lie}{{Lie}}
\DeclareMathOperator{\supp}{{supp}}
\DeclareMathOperator{\Hom}{{Hom}}
\DeclareMathOperator{\Tor}{{Tor}}
\DeclareMathOperator{\Ext}{{Ext}}

\DeclareMathOperator{\Hilb}{{Hilb}}
\DeclareMathOperator{\Poin}{{Poin}}
\DeclareMathOperator{\Prim}{{Prim}}

\DeclareMathOperator{\dist}{{dist}}

\DeclareMathOperator{\Aa}{{A}}
\DeclareMathOperator{\Dd}{{D}}

\newcommand{\surj}{\twoheadrightarrow}

\newcommand{\abs}[1]{\left| #1 \right|}

\newcommand{\ov}[1]{\overline{#1}}

\begin{document}
%% started: September 1, 2004
%% submitted version: December 28, 2004
%% revised version: July 9, 2005

\title[Algebraic invariants for right-angled Artin groups]{%
Algebraic invariants for right-angled Artin groups}

\author[Stefan Papadima]{Stefan Papadima}
\address{Institute of Mathematics of the Academy,
P.O. Box 1-764,
RO-014700 Bucharest, Romania}
\email{Stefan.Papadima@imar.ro}

\author[Alexander~I.~Suciu]{Alexander~I.~Suciu$^1$}
\address{Department of Mathematics,
Northeastern University,
Boston, MA 02115, USA}
\email{a.suciu@neu.edu}
\urladdr{http://www.math.neu.edu/\~{}suciu}

\thanks{$^1$Partially supported by NSF grant DMS-0311142}



\begin{abstract}
A finite simplicial graph $\G$ determines 
a right-angled Artin group $G_\G$, with generators 
corresponding to the vertices of $\G$, and with a relation 
$vw=wv$ for each pair of adjacent vertices.  
We compute the lower central series quotients, the 
Chen quotients, and the (first) resonance variety of 
$G_{\G}$, directly from the graph $\G$.
\end{abstract}
\maketitle

\section{Introduction}
\label{sect:intro}

\subsection{Right-angled Artin groups}
\label{subsec:intro artin}
Let $\G=(\V,\E)$ be a finite simplicial graph. 
To such a graph, there is associated a {\em 
right-angled Artin group}, $G_\G$, with a generator 
$v$ for each vertex $v\in \V$, and with a commutator 
relation $vw=wv$ for each edge $\{v,w\}\in \E$.  

From the graph $\G$, one builds a simplicial complex,  
$\Delta_{\G}$, called the {\em flag complex}, by filling in the 
complete subgraphs on $k$ vertices by $(k-1)$-simplices.  
In turn, the flag complex determines a CW-complex, $K_\G$,
called the {\em cubical complex}, by joining tori in the 
manner prescribed by $\Delta_{\G}$: the $0$-cell corresponds 
to the empty simplex, the $1$-cells correspond to vertices 
$v\in \V$, the $2$-cells correspond to edges $e\in \E$, etc. 

Much is known about the topology of the cubical 
complex $K_\G$, and its fundamental group, $G_\G$.  
First of all, $K_\G$ is an Eilenberg-MacLane $K(G_\G,1)$ 
space, see Charney and Davis \cite{CD} and 
Meier and VanWyk \cite{MV}.  The finiteness properties 
of the kernels of ``diagonal" characters $G_\G \to \Z$ 
have been related in a very precise way to the homotopical 
and homological properties of $\Delta_\G$ by Bestvina 
and Brady \cite{BB}.  Finally, the geometric invariants 
of the type introduced by Bieri, Neumann, and 
Strebel \cite{BNS} (invariants that hold subtle information 
about the finiteness properties of kernels of rational 
characters of $G_\G$) have been determined by 
Meier, Meinert, and VanWyk \cite{MMV}. 

\subsection{Cohomology ring and formality}
\label{subsec:intro coho}
In this paper, we focus on a different set of invariants, mostly 
of an algebraic nature, for the right-angled Artin groups.  
The asphericity of $K_\G$ leads to an identification 
of the cohomology ring $H^*(G_\G,\k)$ with the exterior 
Stanley-Reisner ring $E/J_\G$, where $E$ is the exterior 
algebra over $\k$, on generators $\{ v^*\} _{ v\in \V}$, 
and $J_\G$ is the ideal generated by the monomials 
corresponding to non-edges.  

A crucial observation, due to 
Fr\"oberg \cite{Fr}, is that  $E/J_\G$ is a Koszul algebra.  
When combined with the fact that the group $G_\G$ is 
$1$-formal (cf.~Kapovich and Millson \cite{KM}), this 
implies that the space $K_\G$ is formal, see 
Proposition~\ref{prop:formal kg}. In other words, 
all the rational homotopy theory of $K_\G$ (embodied 
essentially in the Malcev completion of $G_\G$) can be 
extracted, at least in principle, from $E/J_\G$. 

We give here a complete description of several  
invariants of a right-angled Artin group $G_{\G}$---the 
associated graded Lie algebra $\gr(G_{\G})$, the Chen 
groups $\gr_k(G_{\G}/G_{\G}'')$, and the first resonance 
variety $\RR(G_{\G};\k)$---solely in terms of the 
graph $\G$.  Our approach will enable us to derive 
in \cite{PS-bb} similar results for the Bestvina-Brady 
groups associated to $\G$.

\subsection{Lower central series and clique polynomial}
\label{subsec:intro lcs}
We start by analyzing the graded Lie algebra 
$\gr(G)$ associated to a right-angled Artin group $G=G_\G$.  
In Theorem~\ref{thm:lcs artin}, we show that $\gr(G)$ is 
isomorphic, as a graded Lie algebra, with $\HH(G)$, the 
holonomy Lie algebra defined from the reduced diagonal  
map $H_2(G)\to H_1(G)\wedge H_1(G)$.  Furthermore, 
the lower central series quotients are torsion-free, 
with ranks $\phi_k$ given by
\begin{equation}
\label{eq:intro1}
\prod_{k=1}^{\infty}(1-t^k)^{\phi_k}=P_{\G}(-t).  
\end{equation}

Here $P_{\G}(t)=\sum_{k\ge 0} f_k(\G) t^k$ 
is the clique polynomial of $\G$, with coefficient 
$f_k(\G)$ equal to the number of complete $k$-subgraphs 
of $\G$ (in other words, $P_{\G}(t)$ is the $f$-polynomial 
of the flag complex $\Delta_{\G}$).  This LCS formula is a 
manifestation of the Koszul duality between the universal 
enveloping algebra  of $\HH(G)$ and the cohomology ring 
$H^*(G,\k)=E/J$ (this duality was first noted by Shelton 
and Yuzvinsky \cite{SY}).

\subsection{Chen groups and cut polynomial}
\label{subsec:intro chen}
Next, we determine the Chen Lie algebra, $\gr(G/G'')$.  
In Theorem \ref{thm:chen artin} we show that the Chen groups 
of $G=G_\G$ are also torsion-free, with ranks $\theta_k$ given by
\begin{equation}
\label{eq:intro2}
\sum_{k=2}^{\infty} \theta_k t^{k} = Q_{\G}\left( \frac{t}{1-t} \right).  
\end{equation}

Here $Q_{\G}(t)=\sum_{j\ge 2} c_j(\G) t^j$ 
is the ``cut polynomial" of $\G$, with coefficient   
$c_j(\G)$ equal to $\sum_{\W\subset \V\colon  \abs{\W}=j } 
\tilde{b}_0(\G_\W)$, where $\tilde{b}_0(\G_\W)$ is one less 
than the number of components of the full subgraph on $\W$. 

The proof of \eqref{eq:intro2} requires several steps. 
The first step, based on \cite{PS-chen}, is to 
identify, at least rationally, $\gr(G/G'')$ with $\HH(G)/\HH(G)''$.  
The next step, based on \cite{FL} (see also \cite{SS}), 
is to replace the computation of the Chen ranks by 
that of the Betti numbers in the linear strand of the 
free resolution of $E/J$ over $E$.  
A formula from \cite{AAH} and \cite{AHH} now expresses 
these Betti numbers in terms of the corresponding Betti 
numbers of the polynomial Stanley-Reisner ring, $S/I$, 
which in turn can be computed from the cut numbers of $\G$ 
via a well-known formula of Hochster \cite{Ho}.

It is interesting to note that the Chen ranks $\theta_k$ 
often carry more information than the LCS ranks $\phi_k$.  
To illustrate this phenomenon, we give in \S\ref{subsec:circuits} 
examples of graphs  with the same clique polynomial, but 
different cut polynomials.  

\subsection{Resonance varieties and cut sets}
\label{subsec:intro res}

The numerical invariants $\phi_k(G)$ and $\theta_k(G)$ 
only reflect the additive structure of the respective graded 
Lie algebras, $\gr(G)$ and $\gr(G/G'')$.  To capture more 
refined information, related to the Lie bracket, we turn to 
the resonance variety $\RR(G;\k)$.  First defined by Falk \cite{Fa} 
in the context of hyperplane arrangements (over the field $\k=\C$), 
this subvariety of the affine space $H^1(G,\k)$ depends only 
on the holonomy Lie algebra $\HH(G)$.  

It turns out that the resonance variety of $G_{\G}$ is determined 
by the cut sets of $\G$.  More precisely, we show in 
Theorem~\ref{thm:res artin} that 
\begin{equation} 
\label{eq:intro3}
\RR(G_{\G},\k) = \bigcup_{\stackrel{\W\subset \V}{
\G_{\W} \textup{ disconnected}}} H_{\W} , 
\end{equation} 
where $H_{\W}$ is the coordinate subspace 
spanned by $\W$ inside $H^1(G_{\G},\k)= \k^{\abs{\V}}$.  
From this description, it becomes apparent that $\RR(G_{\G},\R)$ 
is the complement of the BNS invariant $\Sigma^1(G_{\G})$.  

To illustrate the strength of resonance, we give in 
\S\ref{subsec:pair} examples of graphs sharing the 
same clique polynomial, and the same cut polynomial, 
but having different resonance varieties.  The difference 
is rather subtle: both varieties have the same number 
of components, of the same dimension, but the respective 
intersection lattices are inequivalent.

\subsection{Higher-dimensional cubical complexes}
\label{subsec:rescaling}

In \cite{KR}, Kim and Roush defined sim\-ply-connected 
analogues of the cubical complex $K_\G$.  The CW-complexes 
$K^q_{\G}$ ($q\ge 1$) are obtained by joining products 
of $(2q+1)$-dimensional spheres in the manner prescribed 
by $\Delta_\G$.  For $q$ sufficiently large, these spaces 
are both formal and coformal. In \S\ref{sect:rescaling}, 
we compute the rational homology and homotopy 
groups of the loop spaces $\Omega K^q_{\G}$, 
purely in terms of the clique polynomial of $\G$. 

We conclude with Example \ref{ex:rescale ad}, which shows once 
again the usefulness of resonance varieties.   Let $\G$ and $\G'$ 
be two trees on the same number of vertices, but with different 
number of extremal vertices.  Then  $K^q_{\G}$ and $K^q_{\G'}$ 
have the same rational homotopy groups, yet non-isomorphic  
homotopy Lie algebras.  The difference in the respective 
Whitehead product structures can be traced back to the 
different number of components in $\RR(G_{\G},\Q)$ and 
$\RR(G_{\G'},\Q)$.

\section{Lie algebras and formality}
\label{sect:lie formal}

\subsection{Holonomy Lie algebra}
\label{subsec:holo lie}
Fix a ground ring $R$, equal to either $\Z$ or a field $\k$.  
Let $A$ be a connected, graded algebra over $R$, with 
graded pieces $A^i$, $i\ge 0$.  Assume $A^1$ is a finitely 
generated $R$-module (torsion-free if $R=\Z$), and 
$A^1\cdot A^1=0$, i.e., $a^2=0$, for all $a\in A^1$.  
Then, the multiplication map in 
degree $1$ descends to an $R$-linear map 
$\mu\colon A^1 \wedge A^1 \to A^{2}$.  

Let $\Lie(A_1)$ be the graded free Lie algebra on  $A_1$, 
where $A_i=(A^i)^{\#}$ denotes the dual $R$-module. 
The {\em holonomy Lie algebra} of $A$, denoted $\HH(A)$, 
is the quotient of $\Lie(A_1)$ by the ideal generated by 
the image of the comultiplication map
\begin{equation}
\label{eq:nabla}
\xymatrix{\nabla\colon A_2 = (A^2)^{\#} \ar[r]^(.35){\mu^{\#}}
&(A^1\wedge A^1)^{\#} \cong A_1\wedge A_1 = 
\Lie_2(A_1)}.
\end{equation}
Note that $\HH(A)$ inherits a natural grading from $\Lie(A_1)$; 
denote by $\HH(A)_k$ the $k$-th graded piece.

Now let $G$ be a finitely presented group.   Then 
$A=H^*(G;R)$ is a connected, graded-commutative 
$R$-algebra, with $A^1$ finitely generated (and torsion-free 
if $R=\Z$).   If $R=\Z$, or $R=\k$ and $\ch\k\ne 2$, then 
automatically $A^1\cdot A^1=0$.  Otherwise, we need to 
assume that the abelianization $H_1(G)$ is torsion-free, 
in which case it can be showed again that 
$A^1\cdot A^1=0$. In either case, we can define 
the holonomy Lie algebra of $G$, with coefficients 
in $R$, to be
\begin{equation}
\label{eq:holo lie}
\HH_R(G)= \HH(H^*(G;R)). 
\end{equation}

In the integral case, simply write $\HH(G)=\HH_{\Z}(G)$, and 
note that $\HH(G)\otimes \Q=\HH_{\Q}(G)$. Moreover, note that, 
in the case when $H_1(G)$ is torsion-free, 
the definition of $\HH(G)$ coincides with 
the one given in \cite{PS-chen}.  

\subsection{Formality properties}
\label{subsec:formal}
In  Appendix A of \cite{Q}, Quillen associates, in a 
functorial way, a Malcev Lie algebra, $M_G$, to any group 
$G$  (see \cite{PS-chen} for further details).  
A finitely presented group $G$ is said to be 
{\em $1$-formal} if $M_G$ is isomorphic to the 
rational holonomy Lie algebra, $\HH_\Q (G)$, 
completed with respect to the bracket length filtration. 
Equivalently, $M_G$ is a quadratic Malcev Lie algebra.

Let $X$ be a finite-type, connected CW-complex.  
Then $X$ is {\em formal}, in the sense of Sullivan \cite{Su}, 
if the rational homotopy type of $X$ is determined 
by the rational cohomology ring $H^*(X;\Q)$. 
As shown by Sullivan, the fundamental group of a 
formal space is $1$-formal.  The $1$-formality of 
$\pi_1(X)$ is not enough to insure the formality of $X$; 
indeed, simply-connected, non-formal spaces are easily 
constructed.  Nevertheless, we have the following 
positive result in this direction. 

Recall that a connected, graded algebra $A$ over 
a field $\k$ is called a {\em Koszul algebra}
if $\Tor_{i}^A(\k,\k)_{j}=0$, for all $i\ne j$.
A necessary condition is that $A$ be the quotient of a
free algebra on generators in degree~$1$ by an
ideal $I$ generated in degree~$2$.  A sufficient 
condition is that $I$ admit a quadratic Gr\"obner basis.

\begin{prop}
\label{prop:formal}
Let $X$ be a finite-type, connected CW-complex.  
Suppose $H^*(X;\Q)$ is a Koszul algebra, and 
$G=\pi_1(X)$ is $1$-formal.  Then $X$ is a formal space.
\end{prop}

\begin{proof}
By assumption, $M_G\cong \widehat{\HH_\Q(G)}$. 
This condition can be reinterpreted in terms of Sullivan's 
minimal models, as follows:
\[
\mathfrak{M}_1(\Omega^*_{\rm dR} X, d) \cong 
\mathfrak{M}_1(H^*(X;\Q),d=0), 
\]
where $\mathfrak{M}(A,d)$ denotes the minimal model 
of a differential graded algebra $(A,d)$, and 
$\mathfrak{M}_1$ denotes the subalgebra generated 
by the degree $1$ part of  $\mathfrak{M}$, see \cite{Su}. 
On the other hand, formality of $X$ means that 
\[
\mathfrak{M}(\Omega^*_{\rm dR} X, d) \cong 
\mathfrak{M}(H^*(X;\Q),d=0).  
\]

The Koszulness assumption on $H^*(X;\Q)$ implies that 
$\mathfrak{M}(\Omega^*_{\rm dR} X, d) = 
\mathfrak{M}_1(\Omega^*_{\rm dR} X, d)$ and 
$\mathfrak{M}(H^*(X;\Q),d=0) =
\mathfrak{M}_1(H^*(X;\Q),d=0)$.  
This observation (which follows readily from 
\cite[Proposition 5.2]{PY}), finishes the proof.
\end{proof}

\subsection{Lower central series }
\label{subsec:lcs}
Let $G$ be a group.  The lower central series of $G$ is the sequence 
of normal subgroups $\{\gamma_k G\}_{k\ge 1}$, defined inductively 
by $\gamma_1 G=G$, $\gamma_2 G=G'$, and 
$\gamma_{k+1}G =(\gamma_k G,G)$, 
where $(x,y)=xyx^{-1}y^{-1}$.  Observe that  
the successive quotients $\gamma_k G/\gamma_{k+1} G$ 
are abelian groups. The direct sum of these quotients, 
\begin{equation}
\label{eq:gr G}
\gr(G)=\bigoplus\nolimits_{k\ge 1} 
\gamma_k G/ \gamma_{k+1} G, 
\end{equation}
is the {\em associated graded Lie algebra} of $G$. 
The Lie bracket $[x,y]$ is induced from the group commutator, 
while the grading is given by bracket length.  

Assume $G$ is finitely generated. Since  
$\gr(G)$ is generated as a Lie algebra by the degree $1$ piece, 
$\gr_1(G)=G/G'$, all the graded pieces 
$\gr_k(G)=\gamma_k G/ \gamma_{k+1} G$ 
are finitely generated abelian groups. 
If, moreover, $G$ is $1$-formal, then
\begin{equation}
\label{eq:holo gr}
 \gr (G)\otimes \Q\cong \HH_\Q(G) 
\end{equation}
as graded Lie algebras, as follows from Appendix A in \cite{Q}.

Assume now $G$ is finitely presented, and $H_1(G)$ is torsion-free.  
Then, the canonical projection $\Lie(H_1(G)) \surj \gr(G)$ 
factors through an epimorphism 
 $\Psi_G \colon  \HH(G) \surj \gr (G)$.  
If, moreover, $G$ is $1$-formal and $\HH(G)$ is torsion-free, 
then the map  $\Psi_G$ yields an isomorphism 
of (integral) graded Lie algebras, 
\begin{equation}
\label{eq:holo gr z}
\xymatrix{\Psi_G\colon \HH(G) \ar[r]^(.55){\cong} & \gr (G)}.
\end{equation}

The Chen Lie algebra of a group $G$ is the associated 
graded Lie algebra of its maximal metabelian quotient, 
$G/G''$.  Assume $G$ is $1$-formal.  Then, as shown 
in \cite{PS-chen}, 
\begin{equation}
\label{eq:holo chen}
\gr (G/G'')\otimes \Q\cong \HH_\Q(G)/\HH_\Q(G)''.
\end{equation}

Assume now $G$ is finitely presented, and $H_1(G)$ is torsion-free.  
Then,  the map $\Psi_G$ descends to an epimorphism 
of graded Lie algebras, 
$\Psi^{(2)}_G \colon  \HH(G)/\HH(G)''  \surj \gr (G/G'')$.  
If, moreover, $G$ is $1$-formal and $ \HH(G)/\HH(G)''$ is torsion-free, 
then, by \cite[Theorem 1.2]{PS-chen}, the map 
$\Psi^{(2)}_G$ yields an isomorphism 
\begin{equation}
\label{eq:holo chen z}
\xymatrix{\Psi_G^{(2)} \colon  \HH(G)/\HH(G)'' 
 \ar[r]^(.6){\cong} & \gr (G/G'')}.
\end{equation}

Of particular interest is the determination of the LCS ranks, 
$\phi_k(G) = \rank \gr_k(G)$ and Chen ranks, 
$\theta_k(G)=\rank \gr_k(G/G'')$.  Note that 
$\phi_1(G)=\theta_1(G)=\rank H_1(G)$; moreover, 
$\phi_2(G)=\theta_2(G)$, $\phi_3(G)=\theta_3(G)$, 
and $\phi_k(G)\ge \theta_k(G)$, for $k\ge 4$. 

 \begin{example}
\label{ex:free group}
Let $F_n$ be the free group of rank $n$.  
The associated graded Lie algebra $\gr(F_n)$ 
is isomorphic to the free Lie algebra $\LL_n=\Lie(\Z^n)$, while 
 $\gr(F_n/F_n'')\cong \LL_n/\LL_n''$.  
The LCS ranks are given by Witt's formula:  
$\phi_k(F_n)=\frac{1}{k}\sum_{d \mid k} \mu(d) n^{k/d}$, 
or, equivalently, 
\begin{equation}
\label{eq:lcs free}
\prod_{k=1}^{\infty}(1-t^k)^{\phi_k(F_n)}=1-n t, 
\end{equation}
while the Chen ranks are given by K.-T. Chen's formula: 
$\theta_k(F_n)=(k-1)\binom{n+k-2}{k}$, for $k\ge 2$, 
or, equivalently, 
\begin{equation}
\label{eq:hilb chen free}
\sum_{k=2}^{\infty}\theta_k(F_n)  t^{k} = 1- 
\frac{1 - n t}{(1-t)^{n}}, 
\end{equation}
\end{example}

\subsection{Homological algebra interpretation}
\label{subsec:homological}
Let $X$ be a finite-type, connected CW complex, 
with fundamental group $G=\pi_1(X)$ and rational 
cohomology ring $A=H^*(X;\Q)$.    
Assume $G$ is a $1$-formal group. 
Using the isomorphism \eqref{eq:holo gr}, Corollary 7.16 from 
Halperin and Stasheff \cite{HS} yields:
\begin{equation}
\label{eq:hs}
\phi_k(G) = \dim_{\Q} \Prim \Ext_A^{k}(\Q,\Q)_k.  
\end{equation}
Here, $\Prim \Ext_A(\Q,\Q)$ denotes the primitives 
in the bigraded Hopf Ext-algebra of $A$, where the 
upper degree is the resolution degree, and the lower degree 
is the internal degree.  

We conclude this section with an analogue of formula \eqref{eq:hs} 
for the Chen ranks.  The key tool is the following result of 
Fr\"oberg and L\"ofwall \cite{FL}.  
   
\begin{prop}
\label{prop:fl}
Let $A$ be a graded algebra over a field $\k$, 
with $A^1$ finite-dimensional and $A^1\cdot A^1=0$. 
Assume $A$ is generated in degree $1$, i.e, the canonical 
map from the exterior algebra $E=\bigwedge A^1$ 
to $A$ is surjective.  Then:
\begin{equation}
\label{eq:fl}
 \left( \HH(A)'/\HH(A)'' \right)_{k}^{\#} = 
\Tor^E_{k-1}(A,\k)_{k}, \quad \text{for $k\ge 2$}.
\end{equation}
\end{prop}

\begin{proof}
In the case when $\ch\k\ne 2$, the conclusion follows 
directly from Theorem 4.1(ii) in \cite{FL}. When $\ch\k=2$, 
the minimal model proof of that Theorem still works, given 
our assumption that $A^1\cdot A^1=0$. 
\end{proof}

\begin{corollary}
\label{cor:chen tor}
Let $G$ be a finitely presented, $1$-formal group.  Let $A$ be a 
graded-commutative $\Q$-algebra, such that $\dim A^1< \infty$ 
and the map $E=\bigwedge A^1 \to A$ is surjective.  
Assume $\HH(A)=\HH_{\Q}(G)$.  Then:
\begin{equation}
\label{eq:chen formula}
\theta_k(G) 
= \dim_{\Q} \Tor^E_{k-1}(A,\Q)_{k},\quad \text{for $k\ge 2$}.
\end{equation}
\end{corollary}

\begin{proof}
Follows from \eqref{eq:holo chen} and Proposition \ref{prop:fl}.
\end{proof}

Given a group $G$ as above, an algebra $A$ with the required 
properties can be constructed as follows.  Let $K=\ker\, (\mu\colon 
H^1(G;\Q)\wedge H^1(G;\Q)\to H^2(G,\Q))$, and $E=\bigwedge 
H^1(G;\Q)$.  Set $A=E/(K+ E^{\ge 3})$.  It is then readily 
checked that  $\HH(A)=\HH_{\Q}(G)$.

In the case when $G$ is the fundamental group of the complement 
of a complex hyperplane arrangement, and $A$ is the Orlik-Solomon 
algebra of the arrangement, a different proof of formula 
\eqref{eq:chen formula} was given in \cite{SS}. 

\section{Homology and lower central series}
\label{sect:artin lcs}

\subsection{The cubical complex}
\label{subsec:cubical complex}
By a graph $\G=(\V_{\G}, \E_{\G})$ we mean a 
loopless, finite graph without multiple edges 
(i.e., a one-dimensional, finite simplicial complex), 
with vertex set $\V_{\G}$ and edge set 
$\E_{\G}\subset \binom{\V_{\G}}{2}$.  To 
such a graph, there corresponds a 
{\em right-angled Artin group} (or, {\em graph group}), 
denoted $G_{\G}$, with presentation 
\begin{equation}
\label{eq:artin group}
G_{\G}= \langle v \in \V_{\G} \mid 
uv=vu \text{ if $\{u,v\} \in \E_{\G}$}\rangle.
\end{equation}

The {\em flag complex} of $\G$, denoted $\Delta_{\G}$, is 
the maximal simplicial complex with $1$-skeleton 
equal to $\G$.  The $k$-simplices of $\Delta_{\G}$ 
are the $(k+1)$-cliques of $\G$, that is, 
the complete subgraphs on $(k+1)$ vertices.  

Let $K_{\G}$ be the associated {\em cubical complex}, 
obtained by joining tori in the manner prescribed by 
$\Delta_{\G}$.  More precisely, for each simplex 
$\sigma$ of $\Delta_{\G}$, let $T_{\sigma}$ be 
the torus formed by identifying parallel faces of a 
$(\abs{\sigma}+1)$-cube ($T_{\emptyset}$ is a point); 
then:
\begin{equation}
\label{eq:cubical}
K_{\G}=\bigcup_{\sigma \in \Delta_{\G}} 
T_{\sigma}\Big\slash \big(T_\sigma\cap T_{\sigma'} =T_{\tau}\  
\text{if}\  \sigma\cap\sigma'=\tau \big)\Big. .
\end{equation}
In other words, if $(S^1)^{\times n}$ is the torus of dimension 
$n=\abs{\V_{\G}}$, with the usual CW-decompo\-sition, 
then $K_{\G}$ is the subcomplex obtained by deleting 
the cells corresponding to the non-faces of $\Delta_{\G}$. 

Note that $H_k(K_{\G})$ is free abelian, of rank 
equal to the number, $f_k(\G)$, of $k$-cliques in $\G$ 
(where $f_0(\G)=1$).   
In particular, the Poincar\'e polynomial of the cubical complex, 
$\Poin(K_{\G},t)=\sum_{k\ge 0} b_k(K_{\G}) t^k$, 
equals the {\em clique polynomial} of the graph, 
$P_{\G}(t)=\sum_{k\ge 0} f_k(\G) t^k$. 

It is readily seen that the fundamental group of $K_{\G}$ 
is the Artin group $G_{\G}$.  In fact, $K_{\G}$ is  
an Eilenberg-MacLane space of type $K(G_{\G},1)$, 
see \cite{CD}, \cite{MV}.

\begin{example}
\label{ex:join}
Let $\G_1$ and $\G_2$ be two graphs (on distinct vertex sets).  
The join of the two graphs, $\G=\G_1* \G_2$, 
is the graph with vertex set $\V=\V_1 \cup \V_2$ and 
edge set $\E=\E_1\cup \E_2 \cup \{ \{v_1,v_2\} \mid v_1\in \V_1 
\text{ and } v_2 \in \V_2\}$.  Clearly, 
$G_{\G}=G_{\G_1}\times G_{\G_2}$ and 
$K_{\G}=K_{\G_1}\times K_{\G_2}$.

For instance, let $\G=K_n$ be the complete graph on $n$ vertices 
(the iterated join of $n$ graphs on a single vertex).  
In this case, $\Delta_{\G}$ is the $(n-1)$-simplex, 
while $K_{\G}$ is the $n$-torus, with fundamental group 
$G_{\G}=\Z^n$. 
\end{example}

\begin{example}
\label{ex:wedge}
Suppose $\G=\G_1\coprod \G_2$ is the disjoint union 
of two subgraphs.  Then 
$G_{\G}=G_{\G_1}* G_{\G_2}$ and 
$K_{\G}=K_{\G_1}\vee K_{\G_2}$.

For instance, let $\G=\overline{K}_n$ be the complement of 
$K_n$ (note that $\overline{K}_n$ is the graph on $n$ vertices, 
with empty edge set).  In this case, $\Delta_{\G}=\overline{K}_n$, 
while $K_{\G}$ is the wedge of $n$ circles, with fundamental 
group $G_{\G}=F_n$. 
\end{example}

For each subset $\W\subset \V$, let 
$\G_{\W}$ be the full subgraph of $\G$ on vertex set $\W$.   
Let $G_{\W}=G_{\G_{\W}}$ be the corresponding 
right-angled Artin group, and let $K_{\W}=K_{\G_{\W}}$ 
be the corresponding cubical complex.  
The inclusion $\W\subset \V$ gives rise to a cellular 
inclusion map $j_{\W}\colon K_{\W}\to K_{\G}$.  
The induced homomorphism, 
$(j_{\W})_{\#} \colon G_{\W}\to G_{\G}$,   
is a split injection, with retract 
 $G_{\G}\to G_{\W}$ given on generators 
 by $v\mapsto v$ if $v\in \W$, and $v\mapsto 1$ otherwise.

\subsection{Cohomology ring and formality}
\label{subsec:coho ring}
Fix a labeling of the vertices of the graph $\G$, and 
write  $\V_{\G}=\{v_1,\dots ,v_n\}$. 
As first noted by Kim and Roush \cite{KR}, the 
cohomology ring of $K_{\G}$ is the quotient of 
$E$, the exterior algebra on generators 
$v_1^*,\dots ,v_n^*$ in degree $1$, 
modulo the ideal $J_{\G}$ generated by the monomials 
$v_i^* v_j^*$ for which $\{v_i,v_j\}$ is not an edge of $\G$.  
In other words, 
\begin{equation}
\label{eq:coho ring}
H^*(K_{\G})=E/J_{\G}
\end{equation}
is the exterior Stanley-Reisner ring of $\G$ (or, 
of the flag complex $\Delta_{\G}$). For example, 
if $\Gamma=K_n$, then $H^*(K_{\G})=E$, while if 
$\Gamma=\overline{K}_n$, then $H^*(K_{\G})=E/\mathfrak{m}^2$, 
where $\mathfrak{m}$ is the ideal generated by $v_1^*,\dots,v_n^*$. 

As noted by Fr\"oberg \cite{Fr} (and reproved by 
Shelton and Yuzvinsky \cite{SY}), 
the set $\{v_i^* v_j^* \mid \{v_i,v_j\}\notin \E_{\G}\}$ is a 
quadratic Gr\"obner basis for $J_{\G}$; thus, 
$H^*(K_{\G},\k)$ is a Koszul algebra, for any 
coefficient field $\k$.  Consequently, $K_{\G}$ 
is not only a $K(\pi,1)$-space, but also 
a rational $K(\pi,1)$-space, i.e., its Bousfield-Kan 
rationalization is aspherical, cf.~\cite{PY}.

As shown by Kapovich and Millson in 
\cite[Theorem 16.10]{KM}, all Artin groups are $1$-formal; 
in particular, the right-angled Artin groups $G_{\G}$ 
are $1$-formal.  Using Proposition \ref{prop:formal}, 
we obtain:

\begin{prop}
\label{prop:formal kg}
Let $K_{\G}$ be the cubical complex associated 
to a graph $\G$. Then $K_{\G}$ is a formal space.
\end{prop}

\subsection{Associated graded Lie algebra}
\label{subsec:gr lie}
We are now ready to compute the graded Lie algebra 
associated to the lower central series filtration of a 
right-angled Artin group. 

\begin{theorem}
\label{thm:lcs artin}
Let $\Gamma=(\V,\E)$ be a finite graph.  
Let $G_\Gamma$ be the corresponding right-angled Artin group, 
with  associated graded Lie algebra $\gr(G_{\Gamma})$.  Then:
\begin{enumerate}
\item \label{i}
The holonomy Lie algebra of $G_{\G}$ has presentation 
\begin{equation}
\label{eq:holo artin}
\HH(G_{\G})=\Lie( \V) / ( [v,w]=0\ \text{if $\{v,w\}\in \E$}).
\end{equation}

\item  \label{ii} 
The graded abelian group $\HH(G_{\G})$ is torsion-free.

\item  \label{iii} 
The canonical projection $\Psi\colon \HH(G_{\G})\to \gr(G_{\G})$ is 
an isomorphism of graded Lie algebras over $\Z$. 

\item  \label{iv}
The ranks of the graded pieces of $\gr(G_{\Gamma})$ are given by
\begin{equation}
\label{eq:lcs artin}
\prod_{k=1}^{\infty}(1-t^k)^{\phi_k(G_{\G})}=
P_{\G}(-t),
\end{equation}
where $P_{\G}(t)=\sum_{k\ge 0} f_k(\G) t^k$ is the clique 
polynomial of the graph $\G$. 
\end{enumerate}
\end{theorem}

\begin{proof}
Part \eqref{i} follows from the definition of the holonomy Lie 
algebra of $G=G_{\G}$, and the identification of $H^*(G)$ 
with the exterior Stanley-Reisner ring of $\G$. 

To prove \eqref{ii}, it is enough to show that the Hilbert series 
of $\HH(G)\otimes \k$ does not depend on the field $\k$.  
By the Poincar\'e-Birkhoff-Witt theorem,  it is enough to 
show that  the Hilbert series of the universal 
enveloping algebra $U \HH_\k(G)$ does not depend 
on $\k$.  Now, as noted by Shelton and Yuzvinsky 
in \cite{SY}, $U \HH_\k(G)= A_{\k}^{!}$, 
where $A_{\k}=H^*(G;\k)$, 
and $A^{!}$ is the Koszul dual of a quadratic algebra $A$. 
But $A_{\k}$ is a Koszul algebra, and so 
$\Hilb(A_{\k}^{!},t)=\Hilb(A_{\k},-t)^{-1}$ 
by Koszul duality.  Finally, note that $\Hilb(A_{\k},t)=P_{\G}(t)$ 
is independent of $\k$. 

Part \eqref{iii} follows from Part \eqref{ii}, together with 
the $1$-formality of $G$ and formula \eqref{eq:holo gr z}. 

Part \eqref{iv} follows from Part \eqref{iii}, together with 
the proof of Part \eqref{ii}.
\end{proof}

\section{Stanley-Reisner rings and Chen groups}
\label{sect:chen artin}

\subsection{Stanley-Reisner rings}
\label{subsec:sr}
Let $\G$ be a finite graph. Put a total order on the vertex set, 
and write $\V_\G=\{v_1,\dots, v_n\}$. As mentioned 
\S\ref{subsec:coho ring}, the cohomology ring $H^*(G_{\G};\k)$ 
is isomorphic to the exterior Stanley-Reisner ring, $E/J_{\G}$, 
where $E=\bigwedge_{\k} (v_1^*,\dots,v_n^*)$ and 
$J_{\G}=\text{ideal} \,\langle v_i^* v_j^* \mid 
\{v_i,v_j\}\notin \E_\G\rangle $.  The Hilbert series of 
this graded ring equals the clique polynomial:  
$\Hilb(E/J_{\G}, t)= P_{\G}(t)$, 

Associated to $\G$ there is also the (polynomial) 
Stanley-Reisner ring, $S/I_{\G}$, where 
$S=\k[x_1,\dots ,x_n]$ and 
$I_{\G}=\text{ideal} \,\langle x_ix_j \mid 
\{v_i,v_j\}\notin \E_\G\rangle $. 
As shown in \cite{Stanley}, the Hilbert series 
of this graded ring is also determined by the clique 
polynomial:
\begin{equation}
\label{eq:hilb sr}
\Hilb(S/I_{\G}, t)= P_{\G}\Big(\frac{t}{1-t}\Big).
\end{equation}

A finitely generated,  $\N^n$-graded module $M$ over 
$A=S$ or $E$ has a 
minimal resolution by free $A$-modules and 
multi-homogeneous $A$-linear maps; such a resolution, 
$F_*\to M$, is unique up to isomorphism. 
Denote by $\beta_{i,a}^A(M)$ the number of basis elements 
in $F_i$ that are homogeneous of multi-degree $a=(a_1,\dots, a_n)$, 
with $a_j\ge 0$.  By minimality of the resolution,
$\beta_{i,a}^A(M)= \dim_\k \Tor^A_{i}(M, \k)_{a}$. 

The multi-graded Poincar\'e series of the two Stanley-Reisner rings 
are related by the following remarkable formula (see \cite[(6.4)]{AHH}, 
\cite[Proposition 2.1]{AAH}, or \cite[Corollary 5.7]{EPY}):
\begin{equation}
\label{eq:ahh-aah}
\sum_{i\ge 0} \sum_{a\in \N^n} \beta_{i,a}^E(E/J_{\G}) t^iu^a = 
\sum_{i\ge 0} \sum_{a\in \N^n} \beta_{i,a}^S(S/I_{\G}) \frac{t^iu^a}
{\prod_{j\in \supp(a)} (1-tu_j)}.
\end{equation}
Moreover, since $I_{\G}$ is a square-free monomial ideal, 
all the nonzero, multi-graded Betti numbers 
$\beta_{i,a}^S(S/I_{\G})$ are also square-free, in the 
sense that $a_j=0$ or $1$; see Hochster \cite{Ho}. 

\subsection{Chen Lie algebra of $G_\G$}
\label{subsec:chen artin}
Before stating our next result, we need to establish some terminology. 
For a graph $\G$ on vertex set $\V$,  let $\kappa(\G)$ be the 
connectivity of $\G$, that is, the maximum integer $r$ so that, 
for any subset $\W \subset \V$ with $\abs{\W}<r$, the full 
subgraph  of $\G$ on $\V\setminus \W$ is connected. 
Also, let $\tilde{b}_0(\G)=\rank \widetilde{H}_0(\G)$ be 
the number of components of $\G$ minus $1$.  

For each $j\ge 1$, define the {\em $j$-th cut number} of $\G$ to be:  
\begin{equation}
\label{eq:cut graph}
c_j(\G)=\sum_{\W\subset \V\colon  \abs{\W}=j } \tilde{b}_0(\G_\W).
\end{equation}
Note that $c_1(\G)=0$, and also $c_j(\G)=0$, if 
$j>  \abs{\V}-\kappa(\G)$. 
We are thus led to define the {\em cut polynomial} of $\G$ as
\begin{equation} 
\label{eq:cut poly}
Q_{\G}(t)=\sum_{j= 2}^{\abs{\V}-\kappa(\G)} c_j(\G) t^j.
\end{equation}

\begin{theorem}
\label{thm:chen artin}
Let $\G$ be a finite graph.  
Let $G_\G$ be the corresponding right-angled Artin group, 
with  Chen Lie algebra $\gr(G_{\G}/G_{\G}'')$.  Then:
\begin{enumerate}
\item  \label{cii} 
The graded abelian group $\HH(G_{\G})/\HH(G_{\G})''$ is torsion-free.

\item  \label{ciii} 
The natural map 
$\Psi^{(2)}\colon \HH(G_\G) / \HH(G_\G)'' \to \gr(G_\G/G_\G'')$ 
is an isomorphism of graded Lie algebras over $\Z$. 

\item  \label{civ}
The ranks of the graded pieces of $\gr(G_{\G}/G_{\G}'')$ 
are given by
\begin{equation}
\label{eq:chen ranks artin}
\sum_{k=2}^{\infty} \theta_k(G_{\G}) t^{k} = 
Q_{\G} \Big(\frac{t}{1-t}\Big),
\end{equation}
where $Q_{\G}(t)=\sum_{j\ge 2} c_j(\G) t^j$ 
is the cut polynomial of $\G$.
\end{enumerate}
\end{theorem}

\begin{proof}
Let $\k$ be a field.  Applying Proposition~\ref{prop:fl} 
to the exterior Stanley-Reisner ring over $\k$, we obtain:
\begin{equation}
\label{eq:frlo}
\dim_\k \left( \HH_\k(G_\G)/\HH_\k(G_\G)'' \right)_k  = 
\beta^{E}_{k-1,k}(E/J_{\G}), \text{ for $k\ge 2$}.  
\end{equation}
Using formula  \eqref{eq:ahh-aah}, we find:
\begin{equation}
\label{eq:aah}
\sum_{k\ge 2} \beta^{E}_{k-1,k}(E/J_{\G})  t^{k} = 
\sum_{i\ge 1} \beta^S_{i,i+1} (S/I_{\G}) \Big(\frac{t}{1-t}\Big)^{i+1}.  
\end{equation}
A well-known formula of Hochster \cite[Theorem 5.1]{Ho} gives:
\begin{equation}
\label{eq:hochster}
\beta^S_{i,i+1} (S/I_{\G}) =\sum_{\W\subset \V\colon \abs{\W}=i+1}
  \dim_\k \widetilde{H}_0 (\G_{\W},\k) = c_{i+1}(\G), 
\end{equation}
compare with \cite[Proposition~2.1]{RV}. 
Combining \eqref{eq:frlo}, \eqref{eq:aah}, and \eqref{eq:hochster}, 
we see that $\dim_\k \left( \HH(G_\G)/\HH(G_\G)'' \right)_k \otimes \k$  
is independent of $\k$. This finishes the proof of Part \eqref{cii}.

Part \eqref{ciii} follows from Part \eqref{cii}, together with 
the $1$-formality of $G_\G$ and formula \eqref{eq:holo chen z}.

Part \eqref{civ}, follows from Part \eqref{ciii}, together with 
formulas \eqref{eq:frlo}--\eqref{eq:hochster}.  
\end{proof}

\section{Cut sets and resonance varieties}
\label{sect:res}

\subsection{Resonance variety of an algebra}
\label{subsec:res var alg}
Let $A$ be a connected, graded algebra over a field $\k$.  
Assume $0<\dim A^1<\infty$, and $A^1\cdot A^1=0$.   
Then, for each $a\in A^1$, right-multiplication by $a$ 
defines a cochain complex
\begin{equation}
\label{A complex}
(A,a)\colon \quad 
\xymatrix{A^0 \ar[r]^{a} & A^1
\ar[r]^{a}  & A^2 }.
\end{equation}

The {\em (first) resonance variety} of $A$ is the locus of 
points $a$ in the affine space $A^1$ where this complex 
fails to be exact in the middle:
\begin{equation} 
\label{eq:resvar}
\RR(A)=\{a \in A^1 \mid H^1(A,a) \ne 0\}.
\end{equation}
Clearly, $\RR(A)$ is a homogeneous algebraic 
variety; in particular, it contains $0\in A^1$.  

Note that a nonzero element $a\in A^1$ 
belongs to $\RR(A)$ if and only if there is an element $a'\in A^1$ 
such that 
\begin{equation} 
\label{eq:res pairs}
a\wedge a' \ne  0 \text{ in } A^1\wedge A^1
\quad \text{and}\quad  \mu(a\wedge a')= 0  \text{ in } A^2,  
\end{equation}
where $\mu\colon A^1 \wedge A^1 \to A^2$ is the 
multiplication map.  In particular, $\RR(A)$ depends 
only on the corestriction of $\mu$ to its image, and hence, 
only on the holonomy Lie algebra $\HH(A)$.  

\begin{lemma}
\label{lem:functorial resonance}
Let $\varphi\colon A\to B$ be a morphism of graded 
algebras.  If $\varphi_1\colon A^1 \to B^1$ is injective, then 
$\varphi_1(\RR(A))\subset  \RR(B)$. 
\end{lemma}

\begin{proof}
We have the following commuting diagram:
\[
\xymatrixcolsep{30pt}
\xymatrix{
A^1 \wedge A^1 \ar[r]^(.6){\mu_A} 
\ar@{^{(}->}[d]^{\varphi_1\wedge \varphi_1} 
& A^2 \ar[d]^{\varphi_2} \\
B^1 \wedge B^1 \ar[r]^(.6){\mu_B} & B^2
}
\]  
Let $a$ be a nonzero element in $\RR(A)\subset A^1$.  
Pick $a'\in A^1$ such that \eqref{eq:res pairs} holds.  Then 
$\varphi_1(a)\wedge \varphi_1(a') \ne 0$, by injectivity 
of $\varphi_1$, and $\mu_B(\varphi_1(a)\wedge \varphi_1(a'))=0$, 
by commutativity of the diagram.  Hence, $\varphi_1(a)$ 
belongs to $\RR(B)$. 
\end{proof}

In particular, if $\varphi\colon A\to B$ is an isomorphism of graded 
algebras, then the linear isomorphism $\varphi_1\colon A^1 \to B^1$ 
restricts to an isomorphism $\RR(A)\to \RR(B)$ between the 
corresponding resonance varieties. 

\begin{lemma}
\label{lem:ortho decomp}
Suppose the vector space $A^1$ decomposes as 
a nontrivial direct sum $U_1\oplus U_2$, with 
$\mu_A(u_1\wedge u_2)=0$ for all $u_1\in U_1$ 
and $u_2\in U_2$.  Then $\RR(A)=A^1$.   
\end{lemma}

\begin{proof}
Let $Z=\{ a\in A^1 \mid a=u_1+u_2, \text { with } 
u_i\in U_i \text{ and } u_i\ne 0\}$.  Clearly, $Z$ is a 
nonempty, Zariski open subset of $A^1$, so it suffices 
to show that $Z\subset \RR(A)$. 

Let $a=u_1+u_2\in Z$.  
Then $a\wedge u_2=u_1\wedge u_2\ne 0$ and 
$\mu_A(a\wedge u_2)=\mu_A(u_1\wedge u_2)=0$.  
Thus, $a$ belongs to $\RR(A)$. 
\end{proof}

\subsection{Resonance variety of a group}
\label{subsec:res var groups}
Let $G$ be a finitely presented group, with $H_1(G)$ nonzero 
and torsion-free. Then, the cohomology ring $A=H^*(G;\k)$ satisfies 
the condition $a^2=0$ for all $a\in A^1$, and so the resonance 
variety of $G$ (over the field $\k$) can be defined as
\begin{equation} 
\label{eq:resvar group}
\RR(G,\k)=\RR(H^*(G;\k)).
\end{equation}

Lemmas \ref{lem:functorial resonance} and \ref{lem:ortho decomp} 
have the following immediate corollaries.  

\begin{corollary}
\label{cor:functorial resonance}
If $\alpha\colon G_1\surj G_2$ is an epimorphism, then 
the induced monomorphism 
$\alpha^*\colon H^1(G_2;\k) \to H^1(G_1;\k)$ takes 
$\RR(G_2,\k)$ to $\RR(G_1,\k)$. 
\end{corollary}

In particular, if $G_1\cong G_2$, there is a linear 
isomorphism of ambient affine spaces, restricting to an 
isomorphism $\RR(G_1,\k)\cong \RR(G_2,\k)$.  
In other words, the ambient linear isomorphism type 
of the resonance variety is an invariant of the group.

\begin{corollary}
\label{cor:ortho decomp}
Suppose $G=G_1* G_2$ splits as a nontrivial free product of 
subgroups $G_1$ and $G_2$, with $H^1(G_i;\k)\ne 0$.  Then 
$\RR(G,\k)=H^1(G;\k)$. 
\end{corollary}

\subsection{Resonance of right-angled Artin groups}
\label{subsec:res artin}

Let $\G=(\V,\E)$ be a finite graph, and let $G_\G$ be the 
corresponding right-angled Artin group. Fix a field 
$\k$. Write $H_{\V}=H^1(G_\G; \k)$ for the $\k$-vector space 
with basis canonically indexed by $\V$.  If $\W$ is a subset of 
$\V$, write $H_{\W}$ for the coordinate subspace spanned 
by $\W$.  

\begin{theorem}
\label{thm:res artin}
Let $\G=(\V,\E)$ be a finite graph.  Then:
\begin{equation} 
\label{eq:res artin}
\RR(G_{\G},\k) = \bigcup_{\stackrel{\W\subset \V}{
\G_{\W} \textup{ disconnected}}} H_{\W}.
\end{equation} 
\end{theorem}

\begin{proof}
Let $U$ be the union of coordinate subspaces on the 
right side of \eqref{eq:res artin}.  We need to show 
$\RR(G_{\G},\k)=U$.  For an element $a=\sum_{v\in \V} a_v v$  
in $H_{\V}$, we will write $\supp(a)=\{v \in \V\mid a_v \ne 0\}$ 
for its support.  

Let $a\in U$. Then there is a subset $\W\subset \V$ 
such that $\supp(a) \subset \W$ and the graph 
 $\G_{\W}$ is disconnected.  (Since $0\in \RR(G_{\G},\k)$, 
 we may assume $a\ne 0$, and so $\supp(a)\ne \emptyset$.)
 Since $\supp(a) \subset \W$, we have $a\in H_{\W}$. 
Since  $\G_{\W}$ is disconnected, $\RR(G_{\W},\k)=H_{\W}$, 
by Corollary \ref{cor:ortho decomp}.  
 Thus, $a\in \RR(G_{\W},\k)$. 
Now,  Corollary \ref{cor:functorial resonance}, 
applied to the natural epimorphism $G_{\G}\surj G_{\W}$, 
shows that $\RR(G_{\W},\k) \subset \RR(G_{\G},\k)$.  
Hence, $a\in \RR(G_{\G},\k)$.  

Let $a\not \in U$.   Then, for every subset $\W\subset \V$ 
containing $\supp(a)$, the graph $\G_{\W}$ is connected. 
We need to show $a\not\in \RR(G_{\G},\k)$, i.e., 
if there is $a'=\sum_{v\in \V} a'_v v \in 
H_{\V}$ such that $\mu(a\wedge a')=0$, then $a\wedge a'=0$.  

Fix a linear order on $\V$. Choose as basis for $H^2(G_{\G},\k)$ 
the set of oriented edges $e=(v,w)$ with $v< w$.  With these 
choices, the multiplication $\mu\colon H_{\V}\wedge H_{\V} \to 
H^2(G_{\G},\k)$ is given by 
\begin{equation}
\label{eq:mu map}
a\wedge a'= \sum_{v<w} 
\left| \begin{smallmatrix} a_v & a_w \\
a'_v & a'_w \end{smallmatrix} \right| v\wedge w
\mapsto 
\mu(a\wedge a')= \sum_{e=(v,w)\in \E} 
\left| \begin{smallmatrix} a_v & a_w \\
a'_v & a'_w \end{smallmatrix} \right|  e.
\end{equation} 
The condition that $\mu(a\wedge a')=0$ is 
thus equivalent to 
\begin{equation}
\label{eq:det}
\begin{vmatrix} a_v & a_w \\
a'_v & a'_w \end{vmatrix}  = 0, 
\quad \text{for all $\{v,w\}\in \E$}.
\end{equation}
To show that $a\wedge a'=0$, it is enough to prove 
the following.

\begin{claim}
Fix $w\in \supp(a)$.  Then 
$a'_v= \frac{a'_{w}}{a_{w}} a_v$, for all $v\in \V$.
\end{claim}

Let $v\in \V$.  First suppose $v\in \supp(a)$.  
By assumption, the full subgraph on $\supp(a)$ is connected.
So pick a path $w=v_0,v_1,\dots, v_k=v$ in $\supp(a)$. 
Using condition \eqref{eq:det} for each edge of this 
path, we find $a'_v=\frac{a'_{v_{k-1}}}{a_{v_{k-1}}} a_v 
= \cdots =  \frac{a'_{v_0}}{a_{v_0}} a_v$, as claimed.

Now suppose $v\not\in \supp(a)$.    
Then, by assumption, the full subgraph 
on $\supp(a) \cup \{v\}$ is connected. 
Hence, $\dist(v, \supp(a))=1$, that is, 
there is a vertex $u\in \supp(a)$ such that 
$\{u, v\}\in \E$.  Condition \eqref{eq:det}  
on this edge implies 
$a_{u} a'_{v} -a_{ v} a'_{u} = 0$.  
Since $a_{v}=0$ and $a_{u}\ne  0$, 
it follows that $a'_{v}=0$, as claimed.
\end{proof}

\begin{corollary}
\label{cor:res var artin}
Let $\G=(\V,\E)$ be a finite graph, and let 
$\RR(G_{\G},\k)\subset H_\V$ be 
the resonance variety of the group $G_\G$.  Then:
\begin{enumerate}
\item 
The irreducible components of 
$\RR(G_{\G},\k)$ are the coordinate subspaces $H_{\W}$, 
maximal among those for which $\G_{\W}$ is disconnected. 

\item 
The codimension of $\RR(G_{\G},\k)$ equals the 
connectivity of $\G$. 
\end{enumerate}
\end{corollary}

In particular, if $\G$ is disconnected, then 
$\RR(G_{\G},\k)=H_{\V}$. 

\begin{remark}
\label{rem:arr groups}
Right-angled Artin groups share some common features 
with arrangement groups, that is, fundamental groups of 
complements of finite unions of hyperplanes in complex 
affine space.  Indeed, if $G$ is a group in 
either class, then $G$ admits a finite presentation, with 
commutator relators only; it is $1$-formal; and  
the resonance variety $\RR(G,\C)$ is a union 
of linear subspaces.  However, the intersection 
of these two classes of groups is fairly small.  
Indeed, a right-angled Artin group $G_{\G}$ is 
an arrangement group if and only if $\G$ is a complete 
multipartite graph.   This assertion will be proved  
elsewhere, in a more general context.  Let us only note 
here that a necessary condition for a group $G$ to be 
an arrangement group is that any two distinct components 
of $\RR(G,\C)$ intersect only at the origin, see \cite{LY}.  
This condition is violated by most graphs we will discuss 
in \S\ref{sect:examples}, for instance, circuits of length 
at least $5$.  
\end{remark}

\subsection{Bieri-Neumann-Strebel invariants}
\label{subsec:bns}

Let $G$ be a finitely generated group. 
Pick a finite generating set for $G$, and let $\mathcal{C}(G)$ 
be the corresponding Cayley graph.  Given an additive real 
character $\chi\colon G\to \R$, let $\mathcal{C}_{\chi}(G)$ 
be the full subgraph on vertex set $\{ g \in G \mid \chi(g)\ge 0\}$.  
In \cite{BNS}, Bieri, Neumann, and Strebel define the (first) 
BNS invariant of $G$ to be:
\begin{equation}
\label{eq:bns inv}
\Sigma^1(G) =\{ \chi \in \Hom(G,\R) \setminus \{0\} \mid 
\text{$\mathcal{C}_{\chi}(G)$ is connected} \}.
\end{equation}
Clearly, $\Sigma^1(G)$ is a conical subset of the vector space 
$\Hom(G,\R)=H^1(G;\R)$.  It turns out that $\Sigma^1(G)$ does 
not depend on the choice of generating set for $G$, see \cite{BNS}.

The BNS invariants of right-angled Artin groups have been 
described in detail by Meier, Meinert, and VanWyk 
\cite{MMV}, \cite{MV}.  Using these descriptions, together 
with our results above, we can identify the set   
$\Sigma^1(G_{\G})$ in terms of the resonance variety 
$\RR(G_{\G},\R)$. 

\begin{prop}
\label{prop:bns and res}
Let $G_{\G}$ be a right-angled Artin group. 
\begin{enumerate}
\item \label{s1}
$\Sigma^1(G_{\G}) = H^1(G_{\G};\R) \setminus \RR(G_{\G},\R)$. 
\item \label{s2}
Suppose $\chi\colon G_{\G}\to \R$ is a rational character, i.e., 
$\im(\chi)\cong \Z$.  Then $\ker (\chi) $ is finitely generated if 
and only if $\chi\notin  \RR(G_{\G},\R)$.
\end{enumerate}
\end{prop}

\begin{proof}
For a character $\chi\colon G_{\G}\to \R$,  let 
$\supp(\chi)=\{v \in \V_\G \mid \chi(v)\ne 0\}$ be its support. 
Theorem \ref{thm:res artin} can be restated, as follows:
$\chi$ is not in $\RR(G_{\G},\R)$ if and only if, for any 
 $\W\subset \V_\G$ with $\supp(\chi)\subset \W$, 
 the graph $\G_{\W}$ is connected. 

In \cite{MV}, Meier and VanWyk show that a character 
$\chi\colon G_{\G}\to \R$ belongs to $\Sigma^1(G_\G)$ if 
and only if the full subgraph on $\supp(\chi)$ is connected 
and dominant.   (A subgraph $\G'\subset \G$ is 
dominant if any vertex $v\in \V_{\G}\setminus \V_{\G'}$ 
has a neighbor in $\G'$.)

Comparing these two characterizations finishes the proof 
of Part \eqref{s1}.

To prove Part \eqref{s2}, we need the following basic 
result from \cite{BNS}, concerning rational characters 
$\chi\colon G\to \R$:  the group $\ker(\chi)$ is 
finitely generated if and only if both $\chi$ and $-\chi$  
belong to $\Sigma^1(G)$.  On the other hand, it 
was shown in \cite{MV} that $\Sigma^1(G_{\G})=-\Sigma^1(G_{\G})$. 
Using these two facts, Part \eqref{s2} follows from Part \eqref{s1}. 
\end{proof}

A generalization to higher $\Sigma$-invariants and higher 
resonance varieties will be given in a forthcoming paper. 

\section{Examples}
\label{sect:examples}

We now illustrate with a few examples how to compute 
the LCS ranks, the Chen ranks, and the resonance varieties 
of right-angled Artin groups, using the results from the 
three previous sections.  We start by discussing methods 
for computing the two graph polynomials that appear in 
formulas \eqref{eq:lcs artin} and \eqref{eq:chen ranks artin}. 

\subsection{Computing the clique and cut polynomials}
\label{subsec:clique cut poly}

When a graph $\G$ has no triangles (i.e., when 
$\Delta_\G=\G$), the clique polynomial is easy 
to compute:   $P_\G(t)=1+ \abs{\V_\G} t + \abs{\E_\G} t^2$.  
For an arbitrary graph, the computation of $P_\G(t)$ is 
more complicated, but a simple recursion formula is 
readily available. 

Fix an edge $e\in \E_{\G}$.  Let $\G\setminus e$ 
be the graph $\G$ with edge $e$ removed, 
and let $\G\setminus N(e)$ be the graph obtained 
from $\G$ by removing the endpoints of $e$, 
all their neighbors, and all edges containing any one of 
these vertices.  Using \eqref{eq:hilb sr} 
and Theorem 5.1 from \cite{Re}, we find:
\begin{equation}
\label{eq:clique recurrence}
P_{\G}(t) = P_{\G \setminus e }(t) - t^2 P_{\G \setminus N(e)}(t).
\end{equation}

Next, we give some (partial) recursion formulas for the  
coefficients of the cut polynomial. 
An edge $e\in \E_{\G}$ is called a {\em bridge} if it does not 
lie in any cycle; we say $e$ is a {\em near-bridge} if 
the only cycles containing $e$ are those incident 
on all vertices of $\G$.

\begin{lemma}
\label{lem:cut deletion}
Suppose the graph $\G$ has a near-bridge $e$. Then, 
for all $j< \abs{\V_{\G}}$, 
\begin{equation}
\label{eq:cut del}
c_j(\G)=c_j(\G \setminus e) -  \binom{\abs{\V_\G}-2}{j-2}.
\end{equation}
\end{lemma}

\begin{proof}
Put $\G'=\G\setminus e$.  Fix a proper subset $\W\subset \V$.  
If $e\not\subset \W$, then clearly $\tilde{b}_0(\G_\W)=\tilde{b}_0(\G'_\W)$.  
If $e\subset \W$, then, by the assumption on $e$, we have 
$\tilde{b}_0(\G_\W)=\tilde{b}_0(\G'_\W)-1$. Hence:
\begin{align*}
\label{eq:cut count}
c_j(\G)&= \sum_{\abs{\W}=j, \, e\not\subset \W} \tilde{b}_0(\G_\W) +
 \sum_{\abs{\W}=j, \, e\subset \W} \tilde{b}_0(\G_\W)\\
 &= \sum_{\abs{\W}=j, \, e\not\subset \W} \tilde{b}_0(\G'_\W) +
 \sum_{\abs{\W}=j, \, e\subset \W}  (\tilde{b}_0(\G'_\W)-1),
\end{align*}
from which \eqref{eq:cut del} follows at once.
\end{proof}

\begin{lemma}
\label{lem:cut single}
Let $\G=\G' \coprod K_1$ be the disjoint union of a graph 
$\G'$ with a singleton graph $K_1$. Then, for all $j\ge 2$:
\begin{equation}
\label{eq:cut single}
c_j(\G)=c_j(\G') + c_{j-1}(\G') + \binom{\abs{\V_\G}-1}{j-1}.
\end{equation}
\end{lemma}

\begin{proof}
We have:
\begin{align*}
\label{eq;cut isolated}
c_j(\G) &= 
\sum_{\abs{\W}=j, \, \W \subset \V_{\G'}} \tilde{b}_0(\G_\W) +
 \sum_{\abs{\W}=j, \, \W \not\subset \V_{\G'}} 
 \tilde{b}_0(\G_\W)\\
 &= 
\sum_{\abs{\W}=j, \, \W \subset \V_{\G'}} \tilde{b}_0(\G'_\W) +
 \sum_{\abs{\W}=j-1, \, \W \subset \V_{\G'}}  (\tilde{b}_0(\G'_\W)+1),
\end{align*}
from which \eqref{eq:cut single} follows at once.
\end{proof}

\subsection{Trees}
\label{subsec:trees}
Let $\G$ be a tree on $n$ vertices.  
In this case, $\Delta_\G=\G$, and so 
the clique polynomial is $P_{\G}(t)=1+nt+(n-1)t^2$.  
Hence, by Theorem \ref{thm:lcs artin}\eqref{iv}, 
\begin{equation}
\label{eq:lcs artin tree}
\prod_{k=1}^{\infty}(1-t^k)^{\phi_k(G_{\G})}=
(1-t)(1-(n-1)t).
\end{equation}

We claim the cut numbers of $\G$ are given by
\begin{equation}
\label{eq:cut tree}
c_{j}(\G)=(j-1) \binom{n-1}{j}, \quad\text{for $j=2,\dots, n-1$}.  
\end{equation}
For $n=1$ or $2$, there is nothing to prove. 
Let $e$ be an extremal edge.  Clearly, $e$ is a bridge, and  
$\G\setminus e = \G_0 \coprod K_1$, where $\G_0$ 
is a tree on $n-1$ vertices.  Formula \eqref{eq:cut tree} 
follows by induction on $n$ (starting from $n=3$, where 
it is obvious), using Lemmas \ref{lem:cut deletion} and 
\ref{lem:cut single}.  

Before proceeding, let us note that \eqref{eq:cut tree} 
is equivalent to $\beta^S_{i,i+1}(S/I_{\G})=
i \binom{n-1}{i+1}$, for $1\le i\le n-2$, 
a fact which can also be deduced from 
\cite[Proposition 11]{Fr85}.

Using Theorem \ref{thm:chen artin}\eqref{civ} and 
formula \eqref{eq:cut tree}, we find that 
the Chen ranks of $G_\G$ are given by 
$\theta_k(G_\G) = (k-1) \cdot \binom{k+n-3}{k}$, 
for all $k\ge 2$, or, equivalently: 
\begin{equation}
\label{eq:hilb chen tree}
\sum_{k=2}^{\infty}\theta_k(G_\G)  t^{k} = 1 - 
\frac{1 - (n-1)t}{(1-t)^{n-1}}.
\end{equation}

Now consider the resonance variety $\RR(G_{\G},\k)$.  The 
minimal cutsets of the tree $\G$ are the non-extremal vertices. 
Thus, $\RR(G_{\G},\k)$ consists of coordinate hyperplanes 
in $\k^n$, one for each non-extremal vertex:
\begin{equation}
\label{eq:res tree}
\RR(G_{\G},\k) = \bigcup_{v\in \V \colon 
v \textup{ non-extremal}} H_{\V\setminus \{v\}} .
\end{equation}

The above discussion shows that, for a tree $\G$, the LCS ranks 
and the Chen ranks of $G_{\G}$ depend only on $n=\abs{V_{\G}}$, 
whereas the resonance variety $\RR(G_{\G},\k)$ also depends on 
the cut sets of the graph.   We illustrate this phenomenon 
with an example. 

\begin{example}
\label{ex:dynkin}
Let $\G$ and $\G'$ be two trees on the same vertex set, but with 
different number of extremal vertices.  (For example, take the 
Dynkin graphs $\G=\Aa_n$ and $\G'=\Dd_n$, with $n\ge 4$; 
clearly, $\G$ has $2$ extremal vertices, whereas $\G'$ has $3$ 
such vertices.)  Then $\phi_k(G_{\G})=\phi_k(G_{\G'})$ and 
$\theta_k(G_{\G})=\theta_k(G_{\G'})$, for all $k\ge 1$, yet 
$\RR(G_{\G},\k)\not\cong \RR(G_{\G'},\k)$.  
Hence,  $\HH(G_{\Aa_n})$ and $\HH(G_{\Dd_n})$ are not isomorphic 
as graded Lie algebras, although they have the same graded ranks. 
\end{example}

\subsection{Circuits}
\label{subsec:circuits}
Let $\G_n$ be a circuit on $n$ vertices, $n\ge 4$.  
Again, $\Delta_{\G_n}=\G_n$, and so $P_{\G_n}(t)=1+nt+nt^2$. 
Hence, 
\begin{equation}
\label{eq:lcs artin circuit}
\prod_{k=1}^{\infty}(1-t^k)^{\phi_k(G_{\G_n})}=
1-nt+nt^2. 
\end{equation}
If $n\ge 5$, the right side of \eqref{eq:lcs artin circuit} 
does not factor into linear factors in $\Z[t]$, and so 
the group $G_{\G_n}$ cannot decompose as an (iterated) semi-direct 
product of free groups, with trivial monodromy action in homology.

Deleting an edge from $\G_n$ yields the $n$-chain $\text{A}_n$. 
Clearly, any edge in $\G_n$ is a near-bridge; thus, by 
Lemma \ref{lem:cut deletion}, 
$c_j(\G_n)=c_j(\Aa_n) -\tbinom{n-2}{j-2}$, for all $j<n$.  
On the other hand, 
$c_n(\G_n)=0$, since $\G_n$ is connected.  
Using formula \eqref{eq:cut tree}, and 
plugging into \eqref{eq:chen ranks artin} gives:
\begin{equation}
\label{eq:hilb chen circuit}
\sum_{k=2}^{\infty}\theta_k(G_{\G_n})  t^{k} = 
\sum_{j=2}^{n-2} \left( (j-1)\binom{n-1}{j} -\binom{n-2}{j-2}\right)
\Big( \frac{t}{1-t}\Big)^{j}.
\end{equation}

The resonance variety $\RR(G_{\G_n}, \k)$ consists of 
$n(n-3)/2$ codimension $2$ subspaces in $\k^n$, 
one for each pair of non-adjacent vertices:
\begin{equation}
 \label{eq:res circuit}
\RR(G_{\G_n},\k) =
\bigcup_{\{v, w\}\in \E_{\overline{\G_n}}} 
H_{\V\setminus \{v,w\}}
\end{equation}

\begin{example}
\label{ex:lcs vs chen}
For each $n\ge 5$, let $\G'_n$ be an $(n-1)$-circuit, 
with an extra vertex and an edge attaching it to the circuit.  
Clearly, $P_{\G'_n}(t)=P_{\G_n}(t)$, and so 
$\phi_k(G_{\G'_n})=\phi_k(G_{\G_n})$, for all $k\ge 1$. 

Although the LCS ranks of $G_{\G'_n}$ and $G_{\G_n}$ 
are the same, the Chen ranks are different.  
Indeed, an application of Lemmas \ref{lem:cut deletion} 
and \ref{lem:cut single} yields $c_j(\G'_n)=c_j(\G_n)$ for 
$j=2,\dots, n-2$.  On the other hand, $c_{n-1}(\G'_n)=1$, 
whereas $c_{n-1}(\G_n)=0$.  Thus, 
\begin{equation}
\label{eq:hilb chen ycircuit}
\sum_{k=2}^{\infty}\theta_k(G_{\G'_n})  t^{k} = 
\sum_{j=2}^{n-2} \left( (j-1)\binom{n-1}{j} -\binom{n-2}{j-2}\right)
\Big( \frac{t}{1-t}\Big)^{j} + \Big( \frac{t}{1-t}\Big)^{n-1}.
\end{equation}
Consequently,  $\theta_k(G_{\G'_n})=\theta_k(G_{\G_n})$ 
for $k\le n-2$, but 
$\theta_k(G_{\G'_n})>\theta_k(G_{\G_n})$, 
for $k\ge n-1$.  
\end{example}

\begin{figure}%
\subfigure{%
\begin{minipage}[t]{0.35\textwidth}
\setlength{\unitlength}{0.6cm}
\begin{picture}(4,4)(-2,-0.6)
\put(3,3){\line(1,-1){3}}
\put(3,3){\line(-1,-1){3}}
\put(0,0){\line(1,0){6}}
\put(1.5,1.5){\line(1,0){3}}
\put(3,0){\line(-1,1){1.5}}
\put(3,0){\line(1,1){1.5}}
\multiput(0,0)(6,0){2}{\circle*{0.3}}
\multiput(1.5,1.5)(3,0){2}{\circle*{0.3}}
\multiput(3,3)(0,-3){2}{\circle*{0.3}}
\put(3,3){\makebox(0,1.1){$1$}}
\put(1.5,1.5){\makebox(-1.1,0){$2$}}
\put(4.5,1.5){\makebox(1.1,0){$3$}}
\put(0,0){\makebox(0,-1.1){$4$}}
\put(3,0){\makebox(0,-1.1){$5$}}
\put(6,0){\makebox(0,-1.1){$6$}}
\end{picture}
\end{minipage}
}
\subfigure{%
\begin{minipage}[t]{0.35\textwidth}
\setlength{\unitlength}{0.6cm}
\begin{picture}(4,4)(-2,-0.6)
\multiput(0,0)(3,0){3}{\line(0,1){3}}
\multiput(0,0)(0,3){2}{\line(1,0){6}}
\put(0,3){\line(1,-1){3}}
\put(3,3){\line(1,-1){3}}
\multiput(0,0)(3,0){3}{\circle*{0.3}}
\multiput(0,3)(3,0){3}{\circle*{0.3}}
\put(0,3){\makebox(0,1.1){$1$}}
\put(3,3){\makebox(0,1.1){$2$}}
\put(6,3){\makebox(0,1.1){$3$}}
\put(0,0){\makebox(0,-1.1){$4$}}
\put(3,0){\makebox(0,-1.1){$5$}}
\put(6,0){\makebox(0,-1.1){$6$}}
\end{picture}
\end{minipage}
}
\caption{\textsf{The graphs $\G$ and $\G'$}}
\label{fig:two triangulations}
\end{figure}

\subsection{A pair of graphs}
\label{subsec:pair}
We conclude this section with a pair of right-angled 
Artin groups that differ in a rather subtle way.

\begin{example}
\label{ex:triangulations}
Let $\G$ and $\G'$ be the two graphs in Figure 
\ref{fig:two triangulations}.  Both graphs have clique polynomial  
$P(t)=1+6t+9t^2+4t^3$ and cut polynomial $Q(t)=t^2(6+8t+3t^2)$; 
thus,  $G_\G$ and $G_{\G'}$  have the same LCS and Chen ranks.  
Their resonance varieties have the same number of irreducible 
components, each a codimension~$2$ linear subspace:
 \[
 \RR(G_{\G},\k)= H_{\ov{23}} \cup  H_{\ov{25}} \cup  H_{\ov{35}}\, , 
 \qquad 
 \RR(G_{\G'},\k)= H_{\ov{15}} \cup  H_{\ov{25}} \cup  H_{\ov{26}}\, .
 \]
Yet the two varieties are not isomorphic, since 
$\dim (H_{\ov{23}} \cap  H_{\ov{25}} \cap  H_{\ov{35}})=3$, 
 whereas $\dim (H_{\ov{15}} \cap  H_{\ov{25}} \cap  H_{\ov{26}})=2$. 
\end{example}

\section{Higher-dimensional cubical complexes}
\label{sect:rescaling}

The cubical complex $K_\G$ associated to a graph $\G$ 
has natural, higher-dimensional analogues, introduced in 
\cite{KR}.  For each integer $q\ge 0$, 
let $K_{\G}^{q}$ be the CW-complex obtained from 
$(S^{2q+1})^{\times n}$, $n=\abs{\V_\G}$, by deleting the cells 
corresponding to the non-faces of the flag complex 
$\Delta_{\G}$. Note that  $K_{\G}^{0}=K_{\G}=K(G_{\G},1)$, 
while $K_{\G}^{q}$ is simply-connected, for $q\ge 1$.  

Let $\Omega Y$ denote the based loop space of a pointed  
CW-complex $Y$. 
In \cite{KR}, Kim and Roush determined the Pontryagin 
ring structure of $H_*(\Omega K_{\G}^{q})$.  We now 
pursue this study, and determine the natural 
Hopf algebra structure on $H_*(\Omega K_{\G}^{q};\Q)$. 
As a byproduct, we compute explicitly the rational homology 
and homotopy groups of $\Omega K_{\G}^{q}$, in terms 
of the clique polynomial of $\G$.

First, some background.  The {\em homotopy Lie algebra} of a 
simply-connected, finite-type CW-complex $Y$ is the graded 
$\Q$-vector space 
\[
\pi_*(\Omega Y) \otimes \Q= 
\bigoplus_{k\ge 1}  \pi_k(\Omega Y)\otimes \Q,
\] 
with the graded Lie algebra structure coming from the 
Whitehead product. By the Milnor-Moore theorem \cite{MM}, 
the universal enveloping algebra of $\pi_*(\Omega Y) \otimes \Q$ 
is isomorphic, as a Hopf algebra, to $H_*(\Omega Y, \Q)$.  
Finally, we say $Y$ is {\em coformal} if its rational homotopy 
type is determined by its homotopy Lie algebra. 

\begin{theorem}
\label{thm:rescale artin}
Let $\G=(\V,\E)$ be a finite graph, and $q$ a positive integer. 
Then:
\begin{enumerate}
\item \label{r1}
The homotopy Lie algebra of $K_{\G}^{q}$ has presentation
\begin{equation*}
\label{eq:homotopy lie}
\pi_*(\Omega K_{\G}^{q}) \otimes \Q = 
\Lie( \V [q] ) / ( [v,w]=0\ \text{if $\{v,w\}\in \E$}),
\end{equation*}
where $\Lie( \V [q] )$ denotes the free Lie algebra on 
generators from $\V$, in degree $2q$. 
\item \label{r2}
The ranks of the graded pieces, $\Phi_k(K_{\G}^{q}) = 
\rank \pi_{k}(\Omega K_{\G}^{q})$, vanish if $2q\nmid k$, 
while the other ranks are given by the following 
homotopy LCS formula:
\begin{equation*}
\label{eq:lcs hlie}
\prod_{k=1}^{\infty}(1-t^{(2q+1)k})^{\Phi_{2qk}(K_{\G}^{q})}=
P_{\G}(-t^{2q+1}).
\end{equation*}
\item \label{r3}
The Hilbert series of $H_*(\Omega K_{\G}^{q}; \Q)$ is given by
\[
\Poin (\Omega K_{\G}^{q},t) = P_{\G}(-t^{2q})^{-1}.
\]
\item \label{r4}
If $2q+1>\dim K_{\G}$, then $K_{\G}^{q}$ is both formal 
and coformal. 
\end{enumerate}
\end{theorem}

\begin{proof}
Let $X=K_{\G}$ and $Y=K^{q}_{\G}$. 
Following \cite{KR}, note that 
the ring $H^*(Y)$ is the $q$-rescaling of the ring 
$H^*(X)$, with rescaling factor of $2q+1$;  
in particular, 
\begin{equation}
\label{eq:xy}
\Poin(Y,t)=\Poin(X,t^{2q+1})=P_{\G}(t^{2q+1}).
\end{equation}

Now recall that $H^*(X;\Q)$ is a Koszul algebra. 
Therefore, by \cite[Theorem~B]{PS-rescale}, 
the graded Lie algebra $\pi_*(\Omega Y)\otimes \Q$ 
is the $q$-rescaling of $\gr(G_{\G})\otimes \Q$, with 
rescaling factor of $2q$.  
Part \eqref{r1} follows from Theorem \ref{thm:lcs artin}, 
Parts \eqref{i} and \eqref{iii}. 

Parts \eqref{r2} and  \eqref{r3} are direct consequences 
of Theorem~C and formula (28) from \cite{PS-rescale}. 

As for Part \eqref{r4}, formality is guaranteed by 
\cite[\S4.1 and Proposition~4.4]{PS-rescale}, 
while coformality follows from \cite[Proposition~1.12]{PS-rescale}.
\end{proof} 

\begin{remark}
\label{rem:non-rescaling}
The conclusion of Theorem  \ref{thm:rescale artin} 
may not hold for an arbitrary 
subcomplex of $(S^{2q+1})^{\times n}$.  Indeed, let $Y$ 
be the $(2q+1)(n-1)$-skeleton of $(S^{2q+1})^{\times n}$. 
Then $Y$ is formal.  Nevertheless, the vanishing property from 
Part \eqref{r2} and the coformality property from Part \eqref{r4} 
fail for $Y$, as soon as $2q+1>n-1>1$; see 
\cite[Example~8.5]{PS-rescale}.
\end{remark}

\begin{example}
\label{ex:rescale ad}
Let $\G$ and $\G'$ be two trees on the same vertex set, but with 
different number of extremal vertices, e.g.,  $\G=\Aa_n$ and 
$\G'=\Dd_n$, with $n\ge 4$.  Let  $Y=K^{q}_{\G}$ and 
$Y'=K^{q}_{\G'}$ be the corresponding higher cubical complexes. 
Then $Y$ and $Y'$ are both formal and coformal spaces, 
by Theorem  \ref{thm:rescale artin}\eqref{r4}, 
and share the same Poincar\'e polynomial, by formula \eqref{eq:xy}; 
moreover, the rational homology and homotopy groups of 
$\Omega Y$ and $\Omega Y'$ have the same dimensions, by 
Theorem  \ref{thm:rescale artin}\eqref{r3} and \eqref{r2}.  
Nevertheless, $Y$ and $Y'$ are not (rationally) 
homotopy equivalent.  Indeed, by the computation from 
Example~\ref{ex:dynkin},  $\HH(G_{\G})\not\cong \HH(G_{\G'})$.  
Comparing the Lie algebra presentations from 
Theorems \ref{thm:lcs artin} and \ref{thm:rescale artin}, 
we conclude  that 
$\pi_*(\Omega Y)\otimes \Q \not\cong \pi_*(\Omega Y')\otimes \Q$. 
\end{example}

\begin{addproof} 
After completing this paper, we became aware of the 
work of Duchamp and Krob \cite{DK1}, \cite{DK2}.  Using 
Theorem~\ref{thm:lcs artin}\eqref{i}, it is readily seen 
that Theorem~2.1 of \cite{DK2} is equivalent to our 
Theorem \ref{thm:lcs artin}\eqref{iii}, while 
Corollary~II.16 and Theorem~III.3 of  \cite{DK1} 
are equivalent to our Theorem~\ref{thm:lcs artin}\eqref{ii} 
and \eqref{iv}, respectively.  Thus, our methods give 
a completely different (and much shorter) proof of 
Duchamp and Krob's results.
\end{addproof}

\begin{ack}
Computations were carried out with the help 
of GAP 4.4 \cite{gap} and Macaulay 2 \cite{GS}.
 
This work was done while the authors were attending the 
program ``Hyperplane Arrangements and Applications'' at the 
Mathematical Sciences Research Institute in Berkeley, California, 
in Fall, 2004.  We thank MSRI for its support and hospitality 
during this stay.  
\end{ack}

\end{document}